\documentclass[11pt]{amsart}

\newtheorem{theo+}              {Theorem}           [section]
\newtheorem{prop+}  [theo+]     {Proposition}
\newtheorem{coro+}  [theo+]     {Corollary}
\newtheorem{lemm+}  [theo+]     {Lemma}
\newtheorem{exam+}  [theo+]     {Example}
\newtheorem{rema+}  [theo+]     {Remark}
\newtheorem{defi+}  [theo+]     {Definition}

\newenvironment{theorem}{\begin{theo+}}{\end{theo+}}
\newenvironment{proposition}{\begin{prop+}}{\end{prop+}}
\newenvironment{corollary}{\begin{coro+}}{\end{coro+}}
\newenvironment{lemma}{\begin{lemm+}}{\end{lemm+}}
\newenvironment{example}{\begin{exam+}}{\end{exam+}}
\newenvironment{remark}{\begin{rema+}}{\end{rema+}}
\newenvironment{definition}{\begin{defi+}}{\end{defi+}}
\theoremstyle{definition}

\newtheorem*{ack}{{\bf Acknowledgments}}

\renewcommand{\Bbb}{\mathbb}

\newcommand{\tr}{\mbox{tr}}
\newcommand{\rank}{\mbox{rank}}

\def\E{/\kern-1.0em \equiv }

\begin{document}
\title[ ]{On the classification of quadratic harmonic morphisms between
Euclidean spaces}
\bigskip
\thanks{*The first author was supported by the Chinese Education  Commission \\
as ``Visiting scholar program (94A04)''.}
\maketitle

\begin{center}
{\bf Ye-Lin Ou*} \\
 Department of Mathematics, Guangxi University for Nationalities, Nanning
530006, P.R.China.\\
Current : Department of Pure Mathematics, School of Mathematics, University of
Leeds, Leeds LS2 9JT, U.K. (until Jan.1996).\\

\bigskip
{\bf John.C.Wood}\\
 Department of Pure Mathematics, School of Mathematics, University of Leeds,
Leeds LS2 9JT, U.K.
\end{center}
\bigskip

\begin{center}
{\bf Abstract}
\end{center}
\bigskip
We give a classification of quadratic harmonic morphisms between Euclidean
spaces (Theorem 2.4)  after proving a Rank Lemma. We also find a correspondence
between umbilical (Definition 2.7) quadratic harmonic morphisms and Clifford
systems. In the case $ {\Bbb R}^{4}\longrightarrow {\Bbb R}^{3} $, we determine
all quadratic harmonic morphisms and show that, up to a constant factor, they
are  all bi-equivalent (Definition 3.2) to the well-known Hopf construction map
and induce harmonic morphisms bi-equivalent to the Hopf fibration ${\Bbb S}^{3}
\longrightarrow {\Bbb S}^{2}$.
\bigskip

{Keywords: Harmonic maps, Harmonic morphisms, Quadratic maps, Clifford
systems.}
\newpage
\section{Quadratic harmonic morphisms and their equations }

\begin{definition}
A map $ \varphi: {\Bbb R}^{m}\longrightarrow {\Bbb R}^{n}$ is called a {\bf
quadratic map} if all of the components of $\varphi$ are quadratic functions (
i.e. homogeneous polynomials of degree 2 ) in $ x_{1},\ldots, x_{m}$. By a {\bf
quadratic harmonic map} (respectively a {\bf quadratic harmonic morphism}) we
mean a harmonic map (respectively a harmonic morphism) which is also a
quadratic map.
\end{definition}
Note that any quadratic harmonic morphism is a non-constant map by our
definition. From the theory of quadratic functions and bilinear forms we know
that a quadratic map $ \varphi: {\Bbb R}^{m}\longrightarrow {\Bbb R}^{n}$ can
always be written as
\begin{equation}
\varphi(X) = ( X^{t}A_{1}X,\ldots, X^{t}A_{n}X)\notag
\end{equation}
where $X$ denotes the column vectors in ${\Bbb R}^{m}$ , $X^{t}$ the transpose
of $ X$ and the $A_{i} \ (i = 1, \ldots, n)$ are symmetric $ m\times m $
matrices (henceforth called { \bf component matrices}).\\

Quadratic harmonic morphisms form a large class of harmonic morphisms between
Euclidean spaces as the following examples show.
\begin{example}{All the following maps are quadratic harmonic morphisms:}\\

$(i)$ {\bf Quadratic harmonic morphisms from orthogonal multiplications.}\\
It is well-known that the standard multiplications $ f : {\Bbb R}^{n} \times
{\Bbb R}^{n} \longrightarrow {\Bbb R}^{n}$, $n = 1,2,4,$ or $8,$ in the real
algebras of real, complex, quaternionic and Cayley numbers are both orthogonal
multiplications and harmonic morphisms. In fact, Baird \cite{Bai83} (Theorem
7.2.7) proves that these are the only possible dimensions for an orthogonal
multiplication $ f : {\Bbb R}^{p} \times {\Bbb R}^{q} \longrightarrow {\Bbb
R}^{n}$ to be a harmonic morphism.\\

$(ii)$ (see \cite{Bai83}) {\bf The Hopf construction maps}  $ F : {\Bbb R}^{n}
\times {\Bbb R}^{n} \longrightarrow {\Bbb R}^{n+1}$. These are defined by\\
\begin{equation}
F(X,Y) =(\left\|X\right\|^{2}-\left\|Y\right\|^{2}, 2f(X,Y)) \notag
\end{equation}
where $f$ is one of the orthogonal multiplications defined in $(i)$.\\

$(iii)$ {\bf Quadratic harmonic morphisms from Clifford systems.}\\
Let $(P_{1},...,P_{n})$ be a Clifford system on ${\Bbb R}^{2m}$, i.e. an
$n$-tuple of symmetric endomorphisms of ${\Bbb R}^{2m}$ satisfying \ $
P_{i}P_{j} + P_{j}P_{i} = 2{\delta}_{ij}I $ for \ $i,j = 1,...n$. Then it
follows from Baird \cite{Bai83}(Theorem 8.4.1) that
\begin{equation}
F(X) = (\langle P_{1}X,X\rangle,\ldots,\langle P_{n}X,X\rangle)\notag
\end{equation}
(where $ \langle , \rangle $ denotes the inner product in Euclidean space)
is a quadratic harmonic morphism with dilation ${\lambda}^{2}(X) =
4\left\|X\right\|^{2}$ for each $X\in{\Bbb R}^{2m}$.\\

$(iv)$ {\bf Quadratic harmonic morphisms from the complete lifts}.\\
Let $ \varphi:{\Bbb R}^{m}\supset U\longrightarrow {\Bbb  R}^{n}$ be a $C^{1}$
map from an open connected subset of ${\Bbb R}^{m} $ into $ {\Bbb  R}^{n} $.
The (real) complete lift (cf.\cite{Ou95A}Definition 2.1) of $\varphi$ is the
map $\overline{\varphi}:{\Bbb R}^{2m}\supset {U\times{\Bbb
R}^{m}}\longrightarrow {\Bbb  R}^{n}$, given by $\overline{\varphi}(X,Y) =
J(\varphi(X))Y$, where $J(\varphi(X))$ is the Jacobian matrix of $\varphi$ at
$X\in U $. It follows from Ou \cite{Ou95A} (Theorem3.3) that {\bf the complete
lift of any quadratic harmonic morphism is again a quadratic harmonic
morphism}.

\end{example}
For some further examples, see Loubeau \cite{Lou95}.
In the rest of this section we will give equations that characterize quadratic
harmonic morphisms between Euclidean spaces.

\begin{lemma}\label{L1}
Let $ \varphi: {\Bbb R}^{m}\longrightarrow {\Bbb R}^{n}$ be a quadratic map
with \\$\varphi(X) = ( X^{t}A_{1}X,\ldots, X^{t}A_{n}X)$. Then $\varphi$ is
harmonic if and only if

\begin{equation}
\tr A_{i} = 0  ,  \;(i = 1,...,n).\label{11}
\end{equation}

\end{lemma}

\begin{proof}
The harmonicity of $\varphi$ is equivalent to the statement that all components
of $\varphi$ are harmonic functions, which is easily seen to be equivalent to
Equation (\ref{11}).
\end{proof}

\begin{proposition}\label{P1}
Let $ \varphi: {\Bbb R}^{m}\longrightarrow {\Bbb R}^{n}$ be a quadratic map
with \\$\varphi(X) = ( X^{t}A_{1}X,\ldots, X^{t}A_{n}X)$. Then $\varphi$ is
horizontally weakly conformal if and only if the following equations hold
\begin{align}
&A_{i}A_{j} + A_{j}A_{i} = 0 , \;(i,j = 1,...n, \ i \neq j),\label{12}\\
&{A_{i}}^{2} = {A_{j}}^{2} , \;(i,j = 1,...n).\label{13}
\end{align}
\end{proposition}

\begin{proof}
For a map $\varphi (X) = ({\varphi}^{1}(X),\ldots,{\varphi}^{n}(X))$ between
Euclidean spaces, horizontal weakly conformality is equivalent to (See
\cite{Fug78},\cite{Ish79A})
\begin{equation}\label{14}
\langle \nabla {\varphi^{i}(X)},\nabla{\varphi^{j}(X)} \rangle = \lambda^{2}(X)
\delta^{ij}
\end{equation}
where $ \delta^{ij}$ is the Kronecker delta and $\nabla\varphi^{i}(X)$ denotes
the gradient of the component function of $\varphi^{i}(X)$.\\

Now for quadratic map $\varphi$, we can calculate its Jacobian matrix as
\begin{equation}\notag
J(\varphi(X)) = \left(\begin{array}{c}
2X^{t}A_{1}\\ \vdots\\ 2X^{t}A_{n}\end{array}\right).
\end{equation}
It is easily seen that Equation (\ref{14}) is equivalent to the following two
equations
\begin{align}
& X^{t}A_{i}A_{j}X \equiv 0 , \; (i,j = 1,...n, i \neq j)\label{15}\\
& X^{t}{A_{i}}^{2}X \equiv X^{t}{A_{j}}^{2}X ,\; (i,j = 1,...n).\label{16}
\end{align}
Since Equations (\ref{15}) and (\ref{16}) are identities of quadratic functions
in $ x_{1},\ldots,x_{m}$, and noting that $A_{i}A_{j}$ is not symmetric in
general we conclude that (\ref{15}) and (\ref{16}) are equivalent to (\ref{12})
and (\ref{13}) respectively. Thus we end the proof of the proposition.
\end{proof}

It is well-known (see \cite{Fug78},\cite{Ish79A}) that a map between Riemannian
manifolds is a harmonic morphism if and only if it is both a harmonic map and a
horizontal weakly conformal map. So by combining Lemma \ref{L1} and Proposition
\ref{P1} we have

\begin{theorem}\label{T1}
A quadratic map $ \varphi: {\Bbb R}^{m}\longrightarrow {\Bbb R}^{n} \ ( m \geq
n)$ with \\
$\varphi(X) = ( X^{t}A_{1}X,\ldots, X^{t}A_{n}X)$ is a harmonic morphism if and
only if
\begin{align}\notag
 &(\ref{11})\; tr A_{i} = 0  ,\; (i = 1,\ldots, n),\\\notag
 &(\ref{12})\; A_{i}A_{j} + A_{j}A_{i} = 0 , \;(i,j = 1,...n,\  i \neq
j),\\\notag
 &(\ref{13})\; {A_{i}}^{2} = {A_{j}}^{2} ,\; (i,j = 1,...n).\notag
\end{align}
\end{theorem}

\section{the classification}

In this section we shall prove the Rank Lemma for quadratic harmonic morphisms
which will be the basis for the classification theorems.

\begin{lemma}{({\bf The Rank Lemma for quadratic harmonic morphisms})}\\
Let $ \varphi: {\Bbb R}^{m}\longrightarrow {\Bbb R}^{n}$ be a quadratic
harmonic morphism with \\$\varphi(X) = ( X^{t}A_{1}X,\ldots, X^{t}A_{n}X)$,
Then\\
(a) All the component matrices $A_{i}$ have the same rank which is an even
number.\\
(b) All the component matrices $A_{i}$ have the same spectrum.
\end{lemma}

\begin{proof}
Suppose that $ \varphi: {\Bbb R}^{m}\longrightarrow {\Bbb R}^{n}$ is a
quadratic harmonic morphism with $\varphi(X) = ( X^{t}A_{1}X,\ldots,
X^{t}A_{n}X)$. Then by Theorem \ref{T1} we have
\begin{equation}
 {A_{i}}^{2} = {A_{j}}^{2} ,\; (i,j = 1,...n)\notag
\end{equation}
which implies that
\begin{equation}
\rank {A_{i}}^{2} = \rank {A_{j}}^{2}. \notag
\end{equation}
The equality of $\rank A_{i}$ now follows from the following\\
{\bf Claim}. For any symmetric matrix $ A $, $\rank {A}^{2} = \rank A. $\\
{\bf Proof of Claim.} It is a standard fact that $A$ can be diagonalized by an
orthogonal matrix $P$, so $P^{-1} A P = D$ is a diagonal matrix. But
\begin{align}\notag
 &\rank {A}^{2} = \rank P^{-1} A^{2} P = \rank P^{-1} A P P^{-1} A P \\\notag
=&\rank {D}^{2} = \rank D = \rank P^{-1} A P = \rank A. \notag
\end{align}
  Now we show that $ \rank A_{i}$ is even. It suffices to do the proof for
quadratic harmonic morphism  ${\varphi}:{\Bbb R}^{m}\longrightarrow {\Bbb
R}^{2}$ with $\varphi (X) =  ( X^{t}A_{1}X, X^{t}A_{2}X )$. After a suitable
choice of orthogonal coordinates, $A_{1}$ assumes the diagonal form
\begin{equation}\label{21}
A_{1} = \left(\begin{array}{ccc}
D_{1} & 0 & 0\\ 0 & -D_{2}& 0\\ 0& 0& 0_{r}
\end{array}\right)
\end{equation}
where $0_{r}$ denotes the $r \times r$ zero matrix, $D_{1} $ is the $k \times
k$  diagonal matrix with entries the positive eigenvalues $ S_{+} = \{ \lambda
_{1},\ldots,\lambda_{k} \}$, and $D_{2} $ is the  $l \times l$ diagonal matrix
with the entries the absolute values of the negative eigenvalues $ S_{-} = \{
\xi_{1},\ldots,\xi_{l} \}$), where $k+l+r = m $.

Using Equations (\ref{12}) and (\ref{13}) we see that $A_{2}$ must have the
form

\begin{equation}\label{22}
A_{2} = \left(\begin{array}{ccc}
0 & B_{1} & 0\\ B_{1}^{t} & 0 & 0\\ 0& 0& 0_{r}
\end{array}\right)
\end{equation}
where $B_{1}$ denotes a $k \times l$ matrix satisfying $D_{1}B_{1} =
B_{1}D_{2}$, which means

\begin{align}\label{23}
\left( \begin{array}{ccc}
\lambda_{1}b_{11} & \ldots & \lambda_{1}b_{1l}\\
\vdots & \ldots & \vdots \\
\lambda_{k}b_{k1} & \ldots & \lambda_{k}b_{kl}
\end{array} \right) =
\left( \begin{array}{ccc}
\xi_{1}b_{11} & \ldots & \xi_{l}b_{1l}\\
\vdots & \ldots & \vdots \\
\xi_{1}b_{k1} & \ldots & \xi_{l}b_{kl}
\end{array} \right).
\end{align}
Since $ \rank A_{1} = \rank A_{2} = k +l $, we see from Equation (\ref{23})
that any $\lambda_{i} \in S_{+}$ must be equal to one of the numbers in $
S_{-}$ otherwise the $i$th row of $B_{1}$ would be zero vector and $\rank A_{2}
< k+l$. This means that $ S_{+}\subset  S_{-}$. A similar reasoning gives $
S_{-}\subset  S_{+}$. Thus we have $ S_{+} = S_{-}$, which means that $ S_{+}$
and $S_{-}$ have equal numbers of the same elements, i.e. $k = l$ and $ S_{+} =
S_{-} = \{ \lambda _{1},\ldots,\lambda_{k} \}$. Thus $ \rank A_{1}= \rank A_{2}
= 2k $ is even, which ends the proof of (a). For (b) we first note, from the
above proof, that the eigenvalues of a component matrix of a quadratic harmonic
morphism must appear in pairs $ \pm \lambda $. On the other hand, it is
elementary that if a symmetric linear transformation $A^{2} $ has an eigenvalue
$\lambda^{2} > 0$ then $A$ must have one of eigenvalues $\pm \lambda$. Now the
rest of the proof follows from the fact that all $A_{i}^!
 {2} \;\;(i = 2,\ldots,n)$ have eigenvalues $\{ \lambda
_{1}^{2},\ldots,\lambda_{k}^{2} \}$, where  $\lambda _{1},\ldots,\lambda_{k}$
are the positive eigenvalues of $A_{1}$.
\end{proof}

\begin{definition}
Let $ \varphi: {\Bbb R}^{m}\longrightarrow {\Bbb R}^{n}$ be a quadratic
harmonic morphism. Then the {\bf $Q$-rank of $\varphi$}, denoted by {\em
$Q$-$\rank (\varphi)$}, is defined to be the rank of its component matrices.
$\varphi$ is said to be {\bf $Q$-nonsingular} if {\em $Q$-$\rank (\varphi) =
m$} , otherwise it is said to be {\bf $Q$-singular}.
\end{definition}

We are now ready to give a characterization of quadratic harmonic morphisms to
$ {\Bbb R}^{2}$.\\

\begin{proposition}\label{P2}
Let $ \varphi: {\Bbb R}^{m}\longrightarrow {\Bbb R}^{2} \; (m \geq 2)$ be a
quadratic harmonic morphism.\\
(i) If $ \varphi $ is $Q$-nonsingular, then $ m = 2k$ for some $ k \in {\Bbb
N}$ and, with respect to suitable orthogonal coordinates in $ {\Bbb R}^{m}$,
$\varphi$ assumes the normal form
\begin{equation}\notag
\varphi (X) = \left(  X^{t} \left( \begin{array}{cc}
D  & 0\\ 0 & -D
\end{array} \right)X,\; X^{t} \left( \begin{array}{cc}
0  & B_{1} \\ B_{1}^{t} & 0
\end{array} \right)X \right)
\end{equation}
where $D$, $B_{1} \in GL({\Bbb R}, k)$, with $D$ diagonal and satisfying
\begin{equation}
\begin{cases}\notag
DB_{1} = B_{1}D\\
B_{1}^{t}B_{1} = D^{2}.
\end{cases}
\end{equation}
(ii) Otherwise $Q$-$\rank(\varphi) = 2k$, for some $ k,\  0 \leq k < m / 2$,
and $ \varphi$ is the composition of an orthogonal projection $\pi : {\Bbb
R}^{m}\longrightarrow {\Bbb R}^{2k}$ followed by a $Q$-nonsingular quadratic
harmonic morphism  ${\varphi}_{1}: {\Bbb R}^{2k}\longrightarrow {\Bbb R}^{2}$.
\end{proposition}

\begin{proof}
Let ${\varphi}:{\Bbb R}^{m}\longrightarrow {\Bbb R}^{2}$ be given by $\varphi
(X) =  ( X^{t}A_{1}X, X^{t}A_{2}X).$ Then from the Rank Lemma we know that
$Q$-rank($\varphi$) is even. If $\varphi $ is $Q$-nonsingular then
$Q$-rank$(\varphi) = m = 2k$. As in the proof of the Rank Lemma, after a
suitable choice of orthogonal coordinates, $A_{1}$ assumes the normal form
\begin{equation}\label{24}
A_{1} = \left(\begin{array}{cc}
D & 0 \\0 & -D
\end{array}\right)
\end{equation}
where $D$ denotes the $k \times k$ diagonal matrix having the positive
eigenvalues of $A_{1}$ as diagonal entries. Then $A_{2}$ must have the form
\begin{equation}\notag
A_{2} = \left(\begin{array}{cc}
0 & B_{1} \\B_{1}^{t} & 0
\end{array}\right)
\end{equation}
with $B_{1} \in GL({\Bbb R}, k)$ satisfying $DB_{1} = B_{1}D$. This, together
with $B_{1}^{t}B_{1} = D^{2}$ given by (\ref{13}) of Theorem \ref{T1}, gives
(i). Now (ii) follows from the fact that if $\varphi $ is $Q$-singular with $
\rank A = 2k < m$, then after a suitable choice of orthogonal coordinates, $A$
takes the form  (\ref{21}) and consequently $B$ the form (\ref{22}).

\end{proof}

Now we give the Classification Theorem for general quadratic harmonic
morphisms.
\begin{theorem}\label{T2}
Let $ \varphi: {\Bbb R}^{m}\longrightarrow {\Bbb R}^{n} \; (m \geq n)$ be a
quadratic harmonic morphism.\\
(I) If $ \varphi $ is $Q$-nonsingular, then $ m = 2k$ for some $ k \in {\Bbb
N}$ and, with respect to suitable orthogonal coordinates in $ {\Bbb R}^{m}$,
$\varphi$ assumes the normal form
\begin{align}\label{25}
\varphi(X) =   & \left( X^{t}\left(
                                 \begin{array}{cc}
                                   D  & 0 \\
                                   0 & -D
                                  \end{array}
                              \right)
                                   X,\; X^{t}
                              \left(
\begin{array}{cc}
                                         0  & B_{1} \\
                                        B_{1}^{t} & 0
                                   \end{array}
                             \right) X ,\ldots,  \right. \\
                &   \left.  X^{t} \left(
                                   \begin{array}{cc}
                                     0  & B_{n-1} \\
                                    B_{n-1}^{t} & 0
                                      \end{array}
                                   \right) X \right).
                     \notag
\end{align}

where $ D, B_{i} \in GL({\Bbb R}, k)$ with $D$ diagonal having the positive
eigenvalues as its diagonal entries satisfy
\begin{equation}\label{26}
\begin{cases}
DB_{i} = B_{i}D\\
B_{i}^{t}B_{i} = D^{2}\\
B_{i}^{t}B_{j} = - B_{j}^{t}B_{i}.\; \; (i,j,= i,...,n-1, i \neq j).
\end{cases}
\end{equation}
(II) Otherwise $Q$-rank$(\varphi) = 2k$ for some $ k,\  0 \leq k < m / 2$, and
$ \varphi$ is the composition of an orthogonal projection $\pi : {\Bbb
R}^{m}\longrightarrow {\Bbb R}^{2k}$ followed by a $Q$-nonsingular quadratic
harmonic morphism  ${\varphi}_{1}: {\Bbb R}^{2k}\longrightarrow {\Bbb R}^{n}$.
\end{theorem}

\begin{proof}
As in the proof of the Rank Lemma, after a suitable choice of orthogonal
coordinates the first component matrix has the form (\ref{21}), and all the
other component matrices $A_{n}$ has the form
\begin{equation}\notag
A_{i+1} = \left(\begin{array}{ccc}
0 & B_{i} & 0\\ B_{i}^{t} & 0 & 0\\ 0& 0& 0_{r}
\end{array}\right),\; i = 1,\ldots, n-1.\notag
\end{equation}
Now if $\varphi$ is $Q$-singular then $Q$-rank$(\varphi) = 2k < m$, for some $
k,\  0 \leq k < m / 2$. It is easily seen that $ \varphi$ is the composition of
an orthogonal projection $\pi : {\Bbb R}^{m}\longrightarrow {\Bbb R}^{2k}$
followed by a $Q$-nonsingular quadratic harmonic morphism  ${\varphi}_{1}:
{\Bbb R}^{2k}\longrightarrow {\Bbb R}^{n}$. Otherwise $\varphi$ is
$Q$-nonsingular in which case $r = 0$. Thus we have the normal form (\ref{25}).
Note that for $ n > 2 $ Equation (\ref{12}) gives the additional Equation
(\ref{26}).
\end{proof}

\begin{corollary}
Any quadratic harmonic morphism is the composition of an orthogonal projection
followed by a $Q$-nonsingular quadratic harmonic morphism from an
even-dimensional space.
\end {corollary}

\begin{remark}
Thus to study quadratic harmonic morphisms it suffices to consider
$Q$-nonsingular ones from even-dimensional spaces.
\end{remark}

\begin{definition}
A quadratic harmonic morphism $ \varphi: {\Bbb R}^{m}\longrightarrow {\Bbb
R}^{n}$ with \\
$\varphi(X) = ( X^{t}A_{1}X,\ldots, X^{t}A_{n}X)$ is said to be {\bf umbilical}
if all the positive eigenvalues of one (and hence all by The Rank Lemma) of its
component matrices are equal.
\end{definition}
There do exist quadratic harmonic morphisms which are not umbilical as the
following example shows.
\begin{example}
It can be checked that ${\varphi}: {\Bbb R}^{8}\longrightarrow {\Bbb R}^{3}$
given by
\begin{align}
\varphi = &( 2x_{1}^{2}+2x_{2}^{2} + 3x_{3}^{2}+3x_{4}^{2}
-2x_{5}^{2}-2x_{6}^{2} -3x_{7}^{2}-3x_{8}^{2},\notag\\
          &4x_{1}x_{5} + 4x_{2}x_{6} + 6x_{3}x_{8} - 6x_{4}x_{7},\notag\\
         -&4x_{1}x_{6}+4x_{2}x_{5}+6x_{3}x_{7} + 6x_{4}x_{8} )\notag
\end{align}
is a quadratic harmonic morphism which is not umbilical since its component
matrices have two distinct positive eigenvalues.
\end{example}
  For more results on constructions of harmonic morphisms into Euclidean spaces
see Ou \cite{Ou95C}.

\section{Quadratic harmonic morphisms and Clifford systems}

\begin{definition}
i) The (n + 1)-tuple $(P_{0},\ldots, P_{n})$ of symmetric endomorphisms of
${\Bbb R}^{2m}$ is called a {\bf Clifford system} on ${\Bbb R}^{2m}$ if
\begin{equation}\notag
 P_{i}P_{j} + P_{j}P_{i} = 2{\delta}_{ij}I \;\; (i,j = 0,1,...,n).
\end{equation}
ii) Let $(P_{0},\ldots, P_{n})$ and $(Q_{0},\ldots, Q_{n})$ be Clifford systems
on ${\Bbb R}^{2p}$ and ${\Bbb R}^{2q}$ respectively, then $(P_{0}\oplus
Q_{0},\ldots, P_{n} \oplus Q_{n})$ is a Clifford system on  ${\Bbb R}^{2p+2q}$,
the so-called {\bf direct sum} of $(P_{0},\ldots, P_{n})$ and $(Q_{0},\ldots,
Q_{n})$ .\\
iii) A Clifford system  $(P_{0},\ldots, P_{n})$ on  ${\Bbb R}^{2m}$ is called
{\bf irreducible } if it is not possible to write  ${\Bbb R}^{2m}$ as a direct
sum of two non-trivial subspaces which are invariant under all $P_{i}$.\\
iv) Two Clifford systems $(P_{0},\ldots, P_{n})$ and $(Q_{0},\ldots, Q_{n})$ on
 ${\Bbb R}^{2m}$ are said to be {\bf algebraically equivalent} if there exists
$ A \in O({\Bbb R}^{2m})$ such that $ Q_{i} = AP_{i}A^{t}$ for all \; $ i =
0,1,\ldots,n $.
\end{definition}

{}From the representation theory of Clifford algebras (see \cite{Hus66}) we
have the following results:\\
{\bf Theorem A.}(See \cite{FerKarMun81})\\
(a) Each Clifford system is algebraically equivalent to a direct sum of
irreducible Clifford systems.\\
(b) An irreducible Clifford  system  $(P_{0},\ldots, P_{n})$ on ${\Bbb R}^{2m}$
exists precisely for the following values of $ n $ and $ m = \delta (n)$:
\newline
\begin{center}
\begin{tabular}{|c|c|c|c|c|c|c|c|c|c|c|}
\hline n & 1 & 2 & 3 & 4 & 5 & 6 & 7 & 8 & \ldots & n+8 \\ \hline
$\delta (n)$ & 1 & 2 & 4 & 4 & 8 & 8 & 8 & 8 & \ldots & 16 $\delta (n)$ \\
\hline
\end{tabular}
\vspace{.25in}
\newline
\end{center}
(c) For $ n \ \E \;0 \; mod \;4$, there exists exactly one algebraically
equivalent class of irreducible Clifford systems. If  $ n \equiv 0 \; mod \;4$,
there are two.\\

\begin{definition}
Let $ \varphi, \tilde{\varphi}: {\Bbb R}^{m}\longrightarrow {\Bbb R}^{n}$ be
two quadratic harmonic morphisms. Then\\
(1) $ \varphi$ and $ \tilde{\varphi}$ are said to be {\bf domain-equivalent} if
there exists an isometry $P$ of $ {\Bbb R}^{m}$ such that $\varphi =
\tilde{\varphi} \circ P$. They are said to be {\bf bi-equivalent} if there
exist isometries  $P$ of $ {\Bbb R}^{m}$ and $G$ of $ {\Bbb R}^{n}$ such that
$\varphi = G^{-1} \circ \tilde{\varphi} \circ P $.\\
(2) The concepts of {\em domain-equivalence} and {\em bi-equivalence} can be
defined similarly for harmonic morphisms between spheres (or, indeed any
Riemannian manifolds).
\end{definition}

Baird has proved (\cite{Bai83} Theorem 8.4.1) that any Clifford system
\\$(P_{0},\ldots, P_{n})$ on ${\Bbb R}^{2m}$ gives rise to a quadratic harmonic
morphism \\$\varphi:{\Bbb R}^{2m}\longrightarrow {\Bbb R}^{n+1}$ defined by
 $$\varphi(X) = (\langle P_{0}X,X \rangle,\;\langle
P_{1}X,X\rangle,\ldots,\langle P_{n}X,X \rangle).$$
It is easy to see that two Clifford systems $(P_{0},\ldots, P_{n})$ and
$(Q_{0},\ldots, Q_{n})$ on  ${\Bbb R}^{2m}$ are algebraically equivalent if and
only if they give rise to domain-equivalent quadratic harmonic morphisms
$\varphi:{\Bbb R}^{2m}\longrightarrow {\Bbb R}^{n+1}$. It is easy to see that
any quadratic harmonic morphism given by a Clifford system is umbilical. We
shall prove that up to a constant factor all umbilical quadratic harmonic
morphisms arise this way.

\begin{theorem}
Up to a homothetic change of coordinates in $ {\Bbb R}^{m}$, any umbilical
quadratic harmonic morphism $ \varphi: {\Bbb R}^{m}\longrightarrow {\Bbb
R}^{n}$ arises from a Clifford system.
\end{theorem}

\begin{proof}
We need only to show that, up to a homothetic change of the coordinates in $
{\Bbb R}^{2k}$, the component matrices of any $Q$-nonsingular umbilical
quadratic harmonic morphism  $ \varphi: {\Bbb R}^{2k}\longrightarrow {\Bbb
R}^{n}$ represent a Clifford system. Indeed it follows from Theorem \ref{T2}
that with respect to suitable orthogonal coordinates in ${\Bbb R}^{2k}$, $
\varphi$ assumes the normal form (\ref{25}) with $D = \lambda Id$, and it is
easily seen that after a change of scale in $ {\Bbb R}^{2k}$ the component
matrices become
\begin{equation}\notag
A_{1} =  \left( \begin{array}{cc}
I_{k}  & 0\\ 0 & -I_{k}
\end{array} \right),\;  A_{2} = \left( \begin{array}{cc}
0  & {\tilde B_{1}} \\ {\tilde B_{1}^{t}} & 0
\end{array} \right)
,\ \ldots,  A_{n} = \left( \begin{array}{cc}
0  & {\tilde B_{n-1}} \\ {\tilde B_{n-1}^{t}} & 0
\end{array} \right)
\end{equation}
with ${\tilde B_{i}} \in O(k)$ satisfying ${\tilde B_{i}}^{t}{\tilde B_{j}} = -
{\tilde B_{j}}^{t} {\tilde B_{i}}.\; \; (i,j,= i,...,n-1, i \neq j)$. It can be
checked that
\begin{equation}\notag
 A_{\alpha}A_{\beta} + A_{\beta}A_{\alpha} = 2{\delta}_{\alpha \beta}I, \;\;
(\alpha,\beta = 1,...,n).
\end{equation}
Which means that the $A_{\alpha}$ represent a Clifford system. This ends the
proof of the theorem.
\end{proof}
\begin{example}
It is easy to check that ${\varphi}: {\Bbb R}^{8}\longrightarrow {\Bbb R}^{5}$
given by
\begin{align}
\varphi(x,y) = ( 3|x|^{2}-3|y|^{2},\;
          &6x_{1}y_{1} - 6x_{2}y_{2} - 6x_{3}y_{3} - 6x_{4}y_{4},\notag\\
          &6x_{1}y_{2} + 6x_{2}y_{1} + 6x_{3}y_{4} - 6x_{4}y_{3},\notag\\
          &6x_{1}y_{3} + 6x_{3}y_{1} + 6x_{4}y_{2} - 6x_{2}y_{4},\notag\\
          &6x_{1}y_{4} + 6x_{4}y_{1} - 6x_{2}y_{3} - 6x_{3}y_{2}). \notag
\end{align}
is an umbilical quadratic harmonic morphism with all positive eigenvalues equal
to $3$. It is also easy to see that it arises from a Clifford system.
\end{example}

In the rest of this section we will determine all quadratic harmonic morphisms
from ${\Bbb R}^{4}$ to ${\Bbb R}^{3}$ and show that they are all bi-equivalent
to some constant multiple $\lambda \varphi_{0}$ of the standard Hopf
construction map and that, up to a change of scale, they all restrict to ${\Bbb
S}^{3} \longrightarrow {\Bbb S}^{2}$ and hence induce bi-equivalent harmonic
morphisms. We thus recover, by simple means, part of a result of Eells and Yiu
\cite{EelYiu94}.
For further results on the existence of quadratic harmonic morphisms see Ou
\cite{Ou95B}
\begin{theorem}
Up to {\bf domain-equivalence}, all quadratic harmonic morphisms $ \varphi:
{\Bbb R}^{4} \longrightarrow {\Bbb R}^{3}$ are of the form
\begin{align}\label{31}
\varphi_{t} = &\lambda ( x_{1}^{2}+x_{2}^{2}-x_{3}^{2}-x_{4}^{2},\\
              &\;2x_{1}x_{3}\cos t + 2x_{1}x_{4}\sin t - 2x_{2}x_{3}\sin t +
2x_{2}x_{4}\cos t,\notag\\
              &2x_{1}x_{3}\sin t - 2x_{1}x_{4}\cos t + 2x_{2}x_{3}\cos t +
2x_{2}x_{4}\sin t )\notag
\end{align}
where $\lambda \neq 0$ and $ t \in[0,2\pi)$. They are all\; {\bf bi-equivalent
} to a constant multiple of the standard Hopf construction map:
\begin{align}\label{32}
\lambda \varphi_{0} = \lambda (
x_{1}^{2}+x_{2}^{2}-x_{3}^{2}-x_{4}^{2},\;&\;2x_{1}x_{3} + 2x_{2}x_{4},
- 2x_{1}x_{4} + 2x_{2}x_{3} ).
\end{align}

\end{theorem}
\begin{proof}
First we note that any quadratic harmonic morphism $ \varphi: {\Bbb R}^{4}
\longrightarrow {\Bbb R}^{3}$ is $Q$-nonsingular since otherwise, $ \varphi$
would be of the form
\begin{equation}
 {\Bbb R}^{4} \overset{\pi}{\longrightarrow} {\Bbb R}^{2}
\overset{\varphi_{1}}{\longrightarrow} {\Bbb R}^{3}\notag
\end{equation}
where $\varphi_{1}$ is a non-constant quadratic harmonic morphism which is
impossible. Next we\\
{\bf Claim:} All quadratic harmonic morphisms $ \varphi: {\Bbb R}^{4}
\longrightarrow {\Bbb R}^{3}$ are umbilical.\\
{\bf Proof of  Claim:} Let $ \varphi:{\Bbb R}^{4} \longrightarrow {\Bbb R}^{3}$
be a quadratic harmonic morphism. Then from Theorem \ref{T2} we have
\begin{align}
D = \left( \begin{array}{cc}
\lambda_{1} & 0\\ 0 & \lambda_{2}
\end{array} \right),\;\;
 B_{1} = \left( \begin{array}{cc}
a_{11} & a_{12}\\ a_{21} & a_{22}
\end{array} \right),\;\;
 B_{2} = \left( \begin{array}{cc}
b_{11} & b_{12}\\ b_{21} & b_{22}
\end{array} \right) \in GL({\Bbb R},2)\notag
\end{align}
satisfying Equation (\ref{26}). Now suppose that $\lambda_{1} \neq
\lambda_{2}$, then by using the first equation of (\ref{26}) we have
\begin{align}
 B_{1} = \left( \begin{array}{cc}
a_{1} & 0\\ 0 & a_{2}
\end{array} \right),\;\;
 B_{2} = \left( \begin{array}{cc}
b_{1} & 0\\ 0 & b_{2}
\end{array} \right) \in GL({\Bbb R},2).\notag
\end{align}
But then the third equation of (\ref{26}) gives $ a_{1}b_{1} = 0 $ and $
a_{2}b_{2} = 0 $, which is impossible since $ B_{1}, B_{2}$ are invertible.
Thus we must have  $\lambda_{1} = \lambda_{2}$.\\
Now any $Q$-nonsingular quadratic harmonic morphism $ \varphi : {\Bbb R}^{4}
\longrightarrow {\Bbb R}^{3}$ can be assumed to be of the form
\begin{align}\label{33}
\varphi = \lambda \left( x_{1}^{2}+x_{2}^{2}-x_{3}^{2}-x_{4}^{2},\;\; X^{t}
\left( \begin{array}{cc}
0  & B_{1} \\ B_{1}^{t} & 0
\end{array} \right)X
,\;\; X^{t} \left( \begin{array}{cc}
0  & B_{2} \\ B_{2}^{t} & 0
\end{array} \right)X \right)
\end{align}
where $ B_{1},\; B_{2}\in O(2) $ satisfy
\begin{equation}\label{34}
B_{1}^{t}B_{2} = - B_{2}^{t}B_{1}.
\end{equation}
without loss of generality, we may assume that $B_{1} =  \left(
\begin{array}{cc}
\cos t  & \sin t \\
-\sin t & \cos t
\end{array}
\right) \in SO(2) ,$ and $B_{2} =  \left( \begin{array}{cc}
\cos \theta  & \sin \theta \\-\sin \theta& \cos \theta
\end{array} \right) \in SO(2) , \ $ or $B_{2} =  \left( \begin{array}{cc}
\sin \theta  & \cos \theta \\ \cos \theta & -\sin \theta
\end{array} \right) \in O(2) \backslash SO(2) $.
It can be checked that for the  second possibility, Equation ~(\ref{34}) has no
solution. whilst for the first possibility Equation~(\ref{34}) is equivalent to
\begin{equation}\label{35}
\begin{cases}
\cos (t-\theta) = - \cos ( \theta - t) \\
\sin (t - \theta) = - \sin (\theta - t)
\end{cases}
\end{equation}
which has solutions $ t - \theta = \theta - t \; \;mod\;2\pi$ i.e.,
\begin{equation}\label{36}
\theta = t - \frac{\pi}{2}\; (mod\;\pi) = t \pm \frac{\pi}{2} \;(mod\;2\pi)
\end{equation}
Inserting (\ref{36}) into (\ref{33}) we have two families
\begin{align*}
(a) \quad \varphi_{t} = \lambda &(x_{1}^{2}+x_{2}^{2}-x_{3}^{2}-x_{4}^{2},\\
&2x_{1}x_{3}\cos t + 2x_{1}x_{4}\sin t - 2x_{2}x_{3}\sin t + 2x_{2}x_{4}\cos
t,\\
&2x_{1}x_{3}\sin t - 2x_{1}x_{4}\cos t + 2x_{2}x_{3}\cos t + 2x_{2}x_{4}\sin t
)
\end{align*}
and
\begin{align*}
(b) \quad \varphi_{t} = \lambda &( x_{1}^{2}+x_{2}^{2}-x_{3}^{2}-x_{4}^{2},\\
&2x_{1}x_{3}\cos t + 2x_{1}x_{4}\sin t - 2x_{2}x_{3}\sin t + 2x_{2}x_{4}\cos
t,\\
-&2x_{1}x_{3}\sin t + 2x_{1}x_{4}\cos t - 2x_{2}x_{3}\cos t - 2x_{2}x_{4}\sin t
).
\end{align*}
However, family (b) can be obtained from family (a) by an orthogonal change of
coordinates in ${\Bbb R}^{3}$. Thus any $Q$-nonsingular quadratic harmonic
morphism $ \varphi : {\Bbb R}^{4} \longrightarrow {\Bbb R}^{3}$ is
domain-equivalent to some $ \varphi_{t}$, whilst $ \varphi_{t} = G^{-1} \circ
\lambda \varphi_{0}$ for $G$ given by
\begin{align}\label{37}
G =  \left( \begin{array}{ccc}
1  & 0 & 0 \\ 0 & \cos t & \sin t\\0 & -\sin t & \cos t
\end{array} \right) \in SO(3).
\end{align}
therefore any $Q$-nonsingular quadratic harmonic morphism $ \varphi : {\Bbb
R}^{4} \longrightarrow {\Bbb R}^{3}$ is bi-equivalent to a multiple of the Hopf
construction map $ \lambda \varphi_{0}$. This ends the proof of the theorem.
\end{proof}

\begin{corollary}
Up to homothety of ${\Bbb R}^{4} $, all quadratic harmonic morphisms $ \varphi
: {\Bbb R}^{4} \longrightarrow {\Bbb R}^{3}$ arise from algebraically
equivalent irreducible Clifford systems on $ {\Bbb R}^{4}$, and they induce
harmonic morphisms ${\Bbb S}^{3} \longrightarrow {\Bbb S}^{2}$ bi-equivalent to
the standard Hopf fibration given by (\ref{32}) with $ \lambda = 1$.
\end{corollary}
\begin{proof}
It is trivial to check that for all $t$, the Clifford systems on  $ {\Bbb
R}^{4}$ represented by
\begin{align*}
&A_{0} = \left(
\begin{array}{cccc}
1 & 0 & 0 & 0 \\
0 & 1 & 0 & 0 \\
0 & 0 & -1 & 0 \\
0 & 0 & 0 & -1
\end{array}
\right)
,\;
A_{1} = \left(
\begin{array}{cccc}
0 & 0 & \cos t & \sin t\\
0 & 0  & -\sin t & \cos t \\
\cos t & -\sin t & 0 & 0 \\
\sin t & \cos t & 0 & 0
\end{array}
\right)
,\notag \\
&A_{2} = \left(
\begin{array}{cccc}
0 & 0 & \sin t & -\cos t \\
0 & 0  & \cos t  & \sin t\\
\sin t & \cos t & 0 & 0 \\
-\cos t & \sin t & 0 & 0
\end{array}
\right),
\end{align*}
are irreducible and are clearly algebraically equivalent. We have seen that,
after a possible change of scale in ${\Bbb R}^{4} $, $ \varphi_{t} = G^{-1}
\circ \varphi_{0}$ for $G$ given by (\ref{37}). Thus for all $t$,\;
$\varphi_{t} $ restricts to ${\Bbb S}^{3} \longrightarrow {\Bbb S}^{2}$ and is
bi-equivalent to the classical Hopf fibration ${\Bbb S}^{3} \longrightarrow
{\Bbb S}^{2}$, which ends the proof of the Corollary.
\end{proof}

\begin{ack} Both authors would like to thank S. Gudmundsson for his comments
which led to significant improvements of the original manuscript. The first
author wishes to thank the  Department of Pure Mathematics, University of
Leeds, U.K., for the hospitality  and generosity he received during a visit
there where this work was done.
\end{ack}

\end{document}